\documentclass[12pt,fleqn]{article}

\usepackage{amsmath}
\usepackage{amsthm}
\usepackage{amssymb}
\usepackage{amsfonts}
\usepackage{graphicx}
\usepackage{hyperref}
\usepackage{color}

\newtheorem{lemma}{Lemma}

\newtheorem{theorem}{Theorem}

\theoremstyle{remark}

\newtheorem{remark}{Remark}
\theoremstyle{definition}

 \numberwithin{equation}{section}

\numberwithin{equation}{section}

\newcounter{comment}
\def\Re {\mathop{\rm Re}\nolimits}

\begin{document}

\title{Application of  uniform asymptotics  to the connection formulas of the fifth Painlev\'{e} equation}

\author{Zhao-Yun Zeng  and Yu-Qiu Zhao\footnote{Corresponding author (Yu-Qiu Zhao).
 {\it{E-mail address:}} {stszyq@mail.sysu.edu.cn} }}
  \date{\hbox{\small {\it{Department of Mathematics, Sun Yat-sen University, GuangZhou
510275, China}}}
}

\maketitle

\begin{abstract}

We apply  the uniform asymptotics method proposed  by Bassom, Clarkson, Law and
McLeod \cite{APC} to a special
Painlev\'{e} V equation, and we provide a simpler and more rigorous proof of the connection formulas
for a special solution of the   equation, which have been established earlier by
McCoy and Tang via the isomonodromy and WKB methods.

\end{abstract}

%%%%%%%%%%%%%%%%%%%%%%%%%%%%%%%%%%%%%%%%%%%%%%%%%%%%%%%%%%%%%%%%%%

\vspace{5mm}

\noindent 2010 \textit{Mathematics Subject Classification}.  33E17, 33C10, 34E05.

\noindent \textit{Keywords and phrases}:  Connection formula; uniform asymptotics; the fifth Painlev\'{e} transcendent; parabolic cylinder function; Bessel function.
%%%%%%%%%%%%%%%%%%%%%%%%%%%%%%%%%%%%%%%%%%%%%%%%%%%%%%%%%%%%%%%%%%%%%%%
\section{Introduction }
We apply and extend the method of uniform asymptotics proposed  by Bassom, Clarkson, Law and
McLeod in \cite{APC}  to a special case of the fifth Painlev\'{e} equation (PV)
\begin{equation}\label{PV}
 \frac{d^2y}{dx^2}=\left (\frac{1}{2y}+\frac{1}{y-1}\right )\left (\frac{dy}{dx}\right )^2-
\frac{1}{x}\frac{dy}{dx}+(1-2\Theta)\frac{y}{x}-\frac{y(y+1)}{2(y-1)},~~\Theta\in\mathbb{R}\setminus \mathbb{Z},
\end{equation}
which appears in the time-dependent correlation functions of the transverse Ising chain at the critical value of the magnetic field
\cite{Mccoy and Perk1, Mccoy and Perk2}. Our main focus will be on the connection formulas of the equation.

The equation \eqref{PV}  is reducible to a special case of the third Painlev\'{e} equation (PIII). Indeed,   if
we put $y(x)=\big(\frac{w(t)-1}{w(t)+1}\big)^2$ with $x=4t$,
then $w(t)$ satisfies  Painlev\'{e} III, namely
\begin{equation}\label{PIII}
 \frac{d^2w}{dt^2}=
\frac{1}{w}\left(\frac{d w}{dt}\right)^2-\frac{1}{t}\frac{d w}{dt}+\frac{1-2\Theta}{t}(w^2-1)+w^3-\frac{1}{w},
\end{equation}
which arose in Ising model studies \cite{Mccoy Tracy}. McCoy, Tracy  and Wu \cite{Mccoy Tracy} derived $\infty\leftrightarrow0$ connection formulas for
a two-parameter class of bounded solutions $w(t;\Theta,\rho)$ of the one-parameter family of Panlev\'{e} III equations \eqref{PIII}.

The special  Painlev\'{e} V equation \eqref{PV}, or the equivalent  Painlev\'{e} III equation \eqref{PIII},  plays
a crucial  role in problems related to random matrices and random processes, orthogonal polynomials, string theory, and in
exactly solvable statistical mechanics and quantum field models. For example,
for the special case $\Theta\rightarrow 0$, if we set $w(t)=-e^{i\phi(t)}$ in \eqref{PIII}, then
$\phi(t)$ satisfies the following equation
\begin{equation}\label{Harmonic equation}
  \big(t\phi'(t)\big)'-2t\sin2\phi(t)+2\sin\phi(t)=0,
\end{equation}
 which appeared  in the problem of classification for rotation surfaces with harmonic inverse mean curvature \cite{Bobe}.
This equation has been  numerically  investigated by Bobenko {\it et al.}\;\cite{Bobe}, and also studied by Andreev and Kitaev
 \cite{FA3}  based on the results obtained in \cite{FA4} by using  isomonodromy deformation and the WKB method.

The problem concerning  one particle density matrix of impenetrable bosons at zero temperature (see Creamer {\it et al.}\;\cite{Creamer}, and Vaidya
and Tracy \cite{Vaida:Tracy})
has been reduced in \cite{Creamer} to a study of  the equation
 \begin{equation}\label{Temperature equ}
   \varphi''(z)=((\varphi'(z))^2-1)\cot\varphi(z)+\frac{1}{z}(1-\varphi'(z)).
 \end{equation}
Creamer {\it et al.}\;\cite{Creamer} have studied the solution of the equation numerically. An analytical study of \eqref{Temperature equ} has carried
out  by Suleimanov
\cite{Sulei}.

It is worth pointing out  that the equation \eqref{Temperature equ}  is equivalent to equation  \eqref{Harmonic equation}  via the transformation
 of the  dependent variables  $e^{i\varphi(z)}=i\tan\frac{\phi(t)}{2}$ and the  independent variables  $z=2it$.

In another special case when $\Theta=\frac{1}{2}$,   we may take the transformation $w(t)=e^{i\psi(t)}$, and then the Painlev\'{e} III
equation (\ref{PIII})
becomes
\begin{equation}\label{sin-Gondor}
  \psi''(t)+\frac{1}{t}\psi'(t)=2\sin2\psi(t).
\end{equation}For a  numerical  study of  the last equation, see Lamb
  \cite{Lamb}, with a connection to the $\pi$-pulse of the sine-Gordon
equation. Flaschka and Newell \cite{Flas:Newell} have also considered equation  \eqref{sin-Gondor}  by using the methods of monodromy preserving
deformation
 and singular integral equations,  and obtained   representation of the one-parameter family of solutions to \eqref{sin-Gondor}  that are holomorphic at
 the origin.

  The same differential equation \eqref{sin-Gondor}, adapting the new variable $s=4t^2$, is a special case
of the following equation
\begin{equation}\nonumber
  \psi''+\frac{1}{s}\psi'=\frac{1}{8s}\sin(2\psi)-\frac{\tilde{\alpha}^2}{4s^2}\frac{\cos\psi}{\sin^3\psi}
\end{equation} with $\tilde{\alpha}=0$,
which is closely related to the asymptotics of Bessel kernel limit of the  Fredholm determinant describing the statistics of the level spacing of the
eigenvalues of Hermitian matrices of large order in a single interval; see Tracy and Widom \cite{Tracy:Widom}.

Now we see that the Painlev\'{e} V  \eqref{PV} is related to various nonlinear equations with statistic physical backgrounds. Yet the main objective of
the present paper, is to calculate the asymptotics of a special solution as  as $x\to 0$ and $x\rightarrow+\infty$, and to justify  connection formulas
between parameters involved in the asymptotic approximations. It is known that there exist solutions to  \eqref{PV}, regular on the positive real line,
with behavior  at the  origin
\begin{equation}\label{Expanding of y at 0}
  y(x)=\frac{i\rho+\Theta}{i\rho-\Theta}\left [1+(1-2\Theta)x+O\left (x^2\right )\right ]\quad\mbox{as}~x\to 0;
\end{equation} cf. \cite{Mccoy and Perk2, MT1},
 where the parameter $\rho$ is not on the imaginary  axis, and such that
 $|{\rm Im}\;\rho|\geqslant\Theta$.

McCoy and  Tang \cite{MT1,MT2} derived $+\infty\leftrightarrow0$ and $\pm i\infty\leftrightarrow0$ connection formulas
 for  two-parameter  solutions of the one-parameter family of Painlev\'{e} V equations  \eqref{PV}, respectively. In this paper, we focus on the
  $+\infty\leftrightarrow0$ connection problem.
 It is known from \cite{Mccoy and Perk2, MT1} that the  solution of  \eqref{PV}  satisfying boundary condition \eqref{Expanding of y at 0} possesses
 the following  asymptotic expansion
 \begin{equation}\label{Expanding of y at infty}
   y(x)=-1+4x^{-\frac{1}{2}}F_1(s)+4x^{-1}\left [2\Theta-1-2F_1^2(s)\right ]+O\left (x^{-\frac{3}{2}}\right )~~\mbox{as}~x \rightarrow+\infty,
 \end{equation}
 where
 \begin{equation} \label{Expression of F}
    F_1(s)=a e^{2is}+b e^{-2is}
 \end{equation}
 with $s=\frac{x}{4}-ab\ln\frac{x}{4}$. Here $a$ and $b$ are independent of $x$, and  $|{\rm Im}\{ab\}|<\frac{1}{4}$.
The following results state that the parameters $a$ and $b$ in  \eqref{Expression of F}  are explicit functions in the parameter $\rho$
in (\ref{Expanding of y at 0}):
\begin{theorem}\label{connection formula theorem}
 \begin{align}
  a(\rho)&=\frac{1}{2\sqrt{\pi}}\Gamma(1+2ic(\rho)){e}^{\frac{i\pi}{4}}{e}^{-6ic(\rho)\ln2}{e}^{\pi c(\rho)}\left [\cos\pi\Theta-\frac{\rho}{\Theta}\sin\pi\Theta\right ],
  \label{Connection formula a}\\
  b(\rho)&=\frac{1}{2\sqrt{\pi}}\Gamma(1-2ic(\rho)){e}^{-\frac{i\pi}{4}}{e}^{6ic(\rho)\ln2}{e}^{\pi c(\rho)}\left [\cos\pi\Theta+\frac{\rho}{\Theta}\sin\pi\Theta\right ],
  \label{Connection formula b}
  \end{align}
  where
  \begin{equation}\label{function in connect}
    c(\rho)=-\frac{1}{4\pi}\ln\left [\left (1+\frac{\rho^2}{\Theta^2}\right )\sin^2\pi\Theta\right ].
  \end{equation}
\end{theorem}
\vskip .3cm

\noindent
\begin{remark}
  From (\ref{Connection formula a}) and (\ref{Connection formula b}) it is readily observed that
 \begin{equation}\label{con ab}
     a(\rho)b(\rho)=c(\rho).
   \end{equation}
  Here use has been made of the fact that $\Gamma(1+2ic)\Gamma(1-2ic)=\frac{4\pi c}{e^{2\pi c}-e^{-2\pi c}}$.
  From \eqref{function in connect} and  (\ref{con ab}), we see that $|{\rm Im}\{ab\}|<\frac{1}{4}$  if and only if  $\rho$ is not on the parts of imaginary axis with
 $|{\rm Im}\rho|\geqslant\Theta$.
 \end{remark}\vskip .3cm

 The above connection formulas have been established by  McCoy and   Tang in
 \cite{MT1} by using the method of isomonodromy deformation and the matching of WKB solutions.
To explain their approach, we  briefly outline the relation of Painlev\'{e} V with  the theory of monodromy preserving deformations of linear
ordinary differential equations with rational coefficients. The reader is referred to \cite{Fokas book, MJ} for more details.

 The Lax pair of the fifth Painlev\'{e} equation, with parameters $\Theta_0=\Theta_1=\Theta$ and $\Theta_{\infty}=0$,
is a system of linear   differential equations for  the  matrix function $\Psi(\lambda, x)$,
\begin{equation}\label{PV:lax pair 1}
\frac{\partial\Psi}{\partial\lambda}=
\left(\begin{array}{cc}
\frac{x}{2}+\frac{v+\frac{\Theta}{2}}{\lambda}-\frac{v+\frac{\Theta}{2}}{\lambda-1}&~~-\frac{u(v+\Theta)}{\lambda}+\frac{yuv}{\lambda-1}\\
\frac{v}{u\lambda}-\frac{v+\Theta}{uy(\lambda-1)}&~~-\frac{x}{2}-\frac{v+\frac{\Theta}{2}}{\lambda}+\frac{v+\frac{\Theta}{2}}{\lambda-1}\end{array}\right)\Psi
\end{equation}
and
\begin{equation}\label{PV:lax pair 2}
  \frac{\partial\Psi}{\partial x}=\left(\begin{array}{cc}
\frac{\lambda}{2}&~~-\frac{u}{x}(v+\Theta-yv)\\\frac{1}{xu}(v-\frac{v+\Theta}{y})&~~-\frac{\lambda}{2}\end{array}\right)\Psi,
\end{equation}
where $y(x)$, $v(x)$ and $u(x)$ satisfy the following isomonodromy deformation system:
\begin{align}
&x\frac{\mathrm{d}y}{\mathrm{d}x}=xy-2v(y-1)^2+2\Theta(y-1),\label{PV: three equations y}\\
&x\frac{{\rm d}v}{{\rm d}x} =yv^2-\frac{1}{y}(v+\Theta)^2,\label{PV: three equations v}\\
&x\frac{{\rm d}}{{\rm d}x}{\rm ln}u =-2v-\Theta+yv+\frac{1}{y}(v+\Theta).
 \label{PV: three equations u}
 \end{align}
Furthermore, $y(x)$ solves
  the fifth Painlev\'{e} equation \eqref{PV}. In this sense,  the Painlev\'{e} V equation is equivalent to the compatibility condition
   $\Psi_{\lambda x}=  \Psi_{x\lambda}$.

The canonical solutions $\Psi_k=\Psi_k(\lambda)$ of  \eqref{PV:lax pair 1}, $k=1,2$,  are determined by
  their asymptotic approximations. In a neighborhood of the irregular singular point $\lambda=\infty$, they
  have the following asymptotic expansion:
 \begin{equation} \label{Psi function expand at infty}
  \Psi_k^{(\infty)}(\lambda)=\left (I+O\left (\frac{1}{\lambda}\right )\right )\exp\left (\frac{\lambda x}{2}\sigma_3
  \right ),~|\lambda|\rightarrow\infty,~x\in \mathbb{R}^+,~k=1,2,
 \end{equation}
respectively  in the Stokes sectors $\Omega_1$ and $\Omega_2$  defined as
\begin{equation*}\Omega_k=\left\{\lambda:-\frac{\pi}{2}+\pi(k-2)<\arg\lambda<\frac{3\pi}{2}+\pi(k-2),
~|\lambda|>R\right\},~~k=1,2, \end{equation*}
with  arbitrary finite positive constant $R$.

 These functions are related by certain Stokes matrices $S_1$ and $S_2$, namely,
 \begin{equation}\label{Connetion  Stokes matrix}
   \Psi_{2}^{(\infty)}(\lambda)=\Psi_{1}^{(\infty)}(\lambda)S_1,~~~~\Psi_{1}^{(\infty)}(\lambda)=\Psi_{2}^{(\infty)}(e^{2\pi i}\lambda)S_2,
 \end{equation}
where
 \begin{equation}\nonumber
    S_1=\left(\begin{array}{cc}
1&~0\\s_1&~1\end{array}\right),~~~~S_2=\left(\begin{array}{cc}
1&~s_2\\0&~1\end{array}\right)
 \end{equation}
and the constants $s_1$ and $s_2$ are called the Stokes multipliers.

The canonical solutions at $\lambda=0$ and $\lambda=1$ are  $2\times2$ unimodular matrices  $\Psi^{(0)}(\lambda,x)$ and $\Psi^{(1)}(\lambda,x)$,
and can be defined as
\begin{equation}\label{Definition of connetion matrix E}
  \Psi_{1}^{(\infty)}(\lambda,x)=\Psi^{(0)}(\lambda,x)E_0,~~\Psi_{1}^{(\infty)}(\lambda,x)=\Psi^{(1)}(\lambda,x)E_1,
\end{equation}
where
$E_0$ and $E_1$ are unimodular constant matrices. Furthermore, the behavior at the singularities are
\begin{align}
\Psi^{(0)}(\lambda,x)&=\tilde{\Psi}^{(0)}(\lambda,x)\lambda^{\frac{1}{2}\Theta\sigma_3}&{\rm as}~\lambda\rightarrow0,
\label{Psi function behavior at 0}\\
\Psi^{(1)}(\lambda,x)&=\tilde{\Psi}^{(1)}(\lambda,x)(\lambda-1)^{\frac{1}{2}\Theta\sigma_3}&{\rm as}~\lambda\rightarrow1,
\label{Psi function behavior at 1}
\end{align}
with  $\tilde{\Psi}^{(0)}(\lambda,x)$  and  $\tilde{\Psi}^{(1)}(\lambda,x)$  being  holomorphic at $\lambda=0$ and $\lambda=1$, respectively.

 The monodromy matrices $E_k$ at the regular singularities $\lambda=k~(k=0,1)$ fulfil  the following cyclic relation:
\begin{equation}\label{cyclic relation}
  E_0^{-1}{e}^{-\pi i\Theta\sigma_3}E_0=S_1S_2E_1^{-1}{e}^{\pi i\Theta\sigma_3}E_1,
\end{equation}
which specifies $E_0$ and $E_1$ up to   left-multiplicative diagonal matrices diag$(d_0,d_0^{-1})$ and diag$(d_1,d_1^{-1})$, respectively.
The connection matrices $E_0$ and $E_1$ are independent of $x$: The isomonodromic condition $\frac{d E_k}{d x}=0$ holds here owing to
the fact that all of $\Psi_{1}^{(\infty)}(\lambda,x)$,  $\Psi^{(0)}(\lambda,x)$ and $\Psi^{(1)}(\lambda,x)$ are  solutions of  \eqref{PV:lax pair 2}.
We note that none of the entries
in matrices $E_0$ and $E_1$ vanishes. In fact, for example, if $(E_0)_{21}=0$, then  from \cite{FA4} we see that
all the Stokes multipliers $s_k=0$, $k=1,2$. From the results in \cite{MJ} for the connection matrices and in view of \eqref{cyclic relation},
we obtain $\rho=i\Theta$,
  which   contradicts   the restriction on $\rho$.

 Appealing  to the isomonodromic
deformation techniques developed in \cite{Its book}, we may use $E_0$ and   $E_1$   to express the parameters $a$ and $b$ in
 \eqref{Expanding of y at infty}   as functions  of the initial parameter $\rho$ in  \eqref{Expanding of y at 0}.
More exactly, the quantities to be used are
\begin{equation}\label{Invariant quantity 1}
  I_0=\frac{(E_0)_{11}(E_0)_{22}}{(E_0)_{12}(E_0)_{21}}
\end{equation}
and
\begin{equation}\label{Invariant quantity 2}
  I_1=\frac{(E_1)_{11}(E_1)_{22}}{(E_1)_{12}(E_1)_{21}}.
\end{equation}

As mentioned earlier, the connection problem  \eqref{Connection formula a}-\eqref{Connection formula b}  was first established by  McCoy and
 Tang \cite{MT1}.
They considered the asymptotic behavior of the solution $\Psi(\lambda, x)$ to the first order system \eqref{PV:lax pair 1}
as $x\rightarrow\infty$ for $\lambda$ being kept away from the turning point $\frac{1}{2}$ and the singular points $0,1$.
Then, they  matched these WKB solutions with the asymptotic approximation at the turning point involving  parabolic cylinder functions,
and the behavior at the singularities  $\lambda=0,1$  involving  the Bessel functions. Eventually  they were capable of  calculating the two invariants
 \eqref{Invariant quantity 1}  and  \eqref{Invariant quantity 2} as $x\rightarrow+\infty$.
The invariants $I_0$ and $I_1$ for small $x$ can be obtained from the connection matrices $E_0$ and $E_1$, which have been done  by Jimbo
 \cite{MJ}; see also \cite{MT1,MT3}, such that
\begin{align}
&I_0=\frac{i\rho+\Theta}{i\rho-\Theta},~~\mbox{and}\label{Invariant quantity 1 small x}\\
&I_1=\frac{i\rho-\Theta}{i\rho+\Theta}.\label{Invariant quantity 2 small x}
\end{align}
The facts that $I_0$ and $I_1$ are independent of $x$ will lead to equalities   involving $a$,  $b$ and the initial parameters. The formulas in
Theorem \ref{connection formula theorem} would then follow.
The procedure, however,  is complicated, and is difficult to
make rigorous; see a comment made in \cite[p.245]{APC}.

In this paper, we shall provide a hopefully simpler and more rigorous derivation of the  formulas
 \eqref{Connection formula a}-\eqref{Connection formula b} by using
the method of uniform asymptotics proposed by Bassom, Clarkson, Law and
McLeod \cite{APC}.
Along the same lines we may find the work of
 Olver
\cite{FO} and Dunster  \cite{Dunster} for coalescing turning points.
Initially in \cite{APC}, PII has been taken as an example to illustrate the method. While the difficulty in extending the techniques for PII to other
transcendents is also acknowledged by the authors of \cite[p.244]{APC}.
Yet the  method has been applied to the connection problems
for  PIII (Sine-Gordon) and PIV; cf Wong and Zhang \cite{WZ1,WZ2}, and has also been used to find the asymptotic
behavior at infinity of the solutions to PIV \cite{Qin1} and PV \cite{Lu1,Lu2}.

We   briefly outline the uniform asymptotics approach  to derive the connection
formulas  \eqref{Connection formula a}-\eqref{Connection formula b}.
First, we will obtain the second-order differential equation   \eqref{sencond order equation 1}  from the Lax pair \eqref{PV:lax pair 1}. Then, we
substitute  the known large-$x$ asymptotic behaviors of $y(x)$ and $v(x)$    into the second-order equation and obtain
an approximate equation  \eqref{sencond order equation 2}. The equation  \eqref{sencond order equation 2}  has
only two coalescing turning points which coalesce with a regular point. Thus, uniform asymptotic solutions are
to be constructed in terms of the parabolic cylinder functions,   uniformly for $\lambda$ on the
Stokes curves as $x\rightarrow\infty$; cf.  Olver \cite{FO}, and see Theorem \ref{Uniform asymptotic theorem} in Section \ref{Uniform infty} below.

Accordingly,  we will calculate the connection matrix via (\ref{Definition of connetion matrix E}) by using  the asymptotics of the fundamental
solutions on the Stokes curves; cf. Theorem \ref{Uniform asymptotic theorem}, respectively as
$|\lambda|\rightarrow\infty$ and $\lambda\rightarrow0$ (or $\lambda\rightarrow 1$),
for $x\rightarrow+\infty$.
 Consequently, we can use   \eqref{Invariant quantity 1}-\eqref{Invariant quantity 2} to obtain $I_0$ and $I_1$
for large $x$; see  \eqref{Invariant quantity 1 large x}-\eqref{Invariant quantity 2 large x} below.
Finally, equating (\ref{Invariant quantity 1 small x}) and (\ref{Invariant quantity 1 large x}) gives $b=b(\rho)$,
and equating (\ref{Invariant quantity 2 small x}) and (\ref{Invariant quantity 2 large x}) gives $a=a(\rho)$.

The difference between the method of uniform asymptotics and the WKB method is that the latter has to match different approximations in different
regions while in the
uniform asymptotics cases  the complicated matching procedure is not needed. The derivation of the method of uniform asymptotics is
rigorous and may lead  to a simpler  argument   with minimal computational efforts.
However, for the uniform asymptotics method, difficulties may arise in   describing  the geometry of the Stokes curves near the  turning points
in the finite plane, such as determining of the correspondence  between $\pm\alpha$ and $\lambda_{1,2}$ in
 \eqref{Define alpha}.

 The rest of the paper is organized as follows. In Section \ref{Uniform infty}, we derive   uniform approximations to the solutions of
 the second-order differential equation obtained from the Lax pair  \eqref{PV:lax pair 1}  as $x\rightarrow+\infty$ by virtue of the
  parabolic cylinder functions on the Stokes curves. The last section is devoted to the
evaluation of  the two invariants $I_0$ and $I_1$  in  \eqref{Invariant quantity 1}-\eqref{Invariant quantity 2} for large $x$.
The proof of Theorem \ref{connection formula theorem} is also provided in that section.

%%%%%%%%%%%%%%%%%%%%%%%%%%%%%%%%%%%%%%%%%%%%%%%%%%%%%%%%%%%%%%%%%%%%%%%%%%%%%%%%%%%%%%%%%

\section{Uniform asymptotics as $x\rightarrow +\infty$ } \label{Uniform infty}
In the present section, we apply the uniform asymptotic method to deal with the large-$x$ behavior of
the second-order differential equation  \eqref{sencond order equation 2}  obtained from the Lax pair  \eqref{PV:lax pair 1}.

First, we   eliminate the function $u(x)$ in  \eqref{PV:lax pair 1}  by taking the following gauge transformation:
\begin{equation}\label{Gauge transformation}
  \hat{\Psi}=M_1u^{-\frac{1}{2}\sigma_3}\Psi.
\end{equation}
As a result, we have
\begin{equation}\label{Lax pair new 1}
  \frac{d\hat{\Psi}}{d\lambda}=\left(\begin{array}{cc}
A&~B\\C&~-A\end{array}\right)\hat{\Psi},
\end{equation}
where
\begin{equation}\label{M1}
  M_1= \frac{1}{\sqrt{2}}\left(\begin{array}{cc}
1&~i\\i&~1\end{array}\right),
\end{equation}
 \begin{equation}\label{Lax pair new 1 entries}
   \begin{split}
    A&=i\left (\frac{v+\frac{\Theta}{2}}{\lambda}-\frac{yv+\frac{1}{y}(v+\Theta)}{2(\lambda-1)}\right),\\
B&=-i\left (\frac{x}{2}+\frac{v+\frac{\Theta}{2}}{\lambda}-\frac{v+\frac{\Theta}{2}}{\lambda-1}\right )
+\frac{1}{2}\left (-\frac{\Theta}{\lambda}+\frac{1}{\lambda-1}\left (yv-\frac{v+\Theta}{y}\right )\right ),\\
C&=i\left (\frac{x}{2}+\frac{v+\frac{\Theta}{2}}{\lambda}-\frac{v+\frac{\Theta}{2}}{\lambda-1}\right )
+\frac{1}{2}\left (-\frac{\Theta}{\lambda}+\frac{1}{\lambda-1}(yv-\frac{v+\Theta}{y})\right ).
   \end{split}
 \end{equation}

It is easily verified  that the gauge transformation  \eqref{Gauge transformation}  does not change the monodromy matrix $E_{k}$, $k=0,1$.

Let $(\psi_1,\psi_2)^T$ be a vector  solution of  \eqref{Lax pair new 1}, and
set
\begin{equation}\label{Transform psi to phi}
  \phi=B^{-\frac{1}{2}}\psi_1,
\end{equation}
then $\phi(\lambda)$ solves the equation
\begin{equation}\label{sencond order equation 1}
  \frac{d^2\phi}{d\lambda^2}=\left [A^2+BC+A'-AB^{-1}B'+\frac{3}{4}(B^{-1}B')^2-\frac{1}{2}B^{-1}B''\right ]\phi,
\end{equation}where the derivatives are taken with respect to $\lambda$, for instance, $A'=\frac{d A}{d\lambda}$.

 Combining  \eqref{PV: three equations y}  and  \eqref{PV: three equations v}  with \eqref{Expanding of y at infty}, we have
 \begin{equation}\label{Expanding of v at infty}
   v(x)=-\frac{x}{8}-\frac{i}{4}x^{\frac{1}{2}}F_2(s)+\frac{1}{2}\left [F_1^2(s)-\Theta\right ]+O\left (x^{-\frac{1}{2}}\right )
   ~~\mbox{as}~x\rightarrow+\infty,
  \end{equation}
where $F_1(s)$ is defined in  \eqref{Expression of F}, $s=\frac{x}{4}-ab\ln\frac{x}{4}$, and
\begin{equation}\label{Express F2}
  F_2(s)=a e^{2is}-b e^{-2is}.
\end{equation}
For $x\rightarrow+\infty$, substituting  \eqref{Lax pair new 1 entries}  into \eqref{sencond order equation 1},
and  in view of  \eqref{Expanding of y at infty}  and  \eqref{Expanding of v at infty},   a  straightforward calculation yields the following second-order
 equation:
\begin{equation}\label{sencond order equation 2}
  \frac{d^2\phi}{d\lambda^2}=-x^2\left \{-\frac{(2\lambda-1)^2}{16\lambda(\lambda-1)}+\frac{Q_1(\lambda,x)}{x}
-\frac{Q_2(\lambda,x)}{x^{3/2}}+O\left (x^{-2}\right )\right \}\phi:=-x^2Q(\lambda,x)\phi,
\end{equation}
where
\begin{align}
Q_1(\lambda,x)&=\frac{1}{4\lambda(\lambda-1)}\left [4ab+i-i\left (\lambda-\frac{1}{2}\right )\frac{ H' }{H}\right],\label{Q1 function}\\
Q_2(\lambda,x)&=\frac{1}{2\lambda(\lambda-1)}\left [F_1(s)+F_2(s)-\left (\lambda-\frac{1}{2}\right )F_2(s)\frac{ H' }{H}\right ],\label{Q2 function}
\end{align}
with
\begin{equation}\label{expression of H}
  H(\lambda,x)=ix^{1/2}\left (\lambda-\frac{1}{2}\right )^2+
\lambda F_1(s)-\frac{1}{2}F_2(s)+O(x^{-1/2})~~{\rm as}~x\rightarrow\infty.
\end{equation}
Here $F_1(s)$ and $F_2(s)$ are given in \eqref{Expression of F}  and  \eqref{Express F2}, respectively.

For large $x$, equation \eqref{sencond order equation 2} has two turning points
\begin{equation}\label{Two turning points}
  \lambda_j=\frac{1}{2}\mp x^{-1/2}\sqrt{4ab+i}\; (1+o(1)),~~j=1,2,
\end{equation}
  which  coalesce with each other at $\frac{1}{2}$ as $x\rightarrow\infty$.  The Stokes  curves are defined as
\begin{equation}\label{Stokes curves}
  \Re \sqrt{\lambda(\lambda-1)} =0.
\end{equation}

In line with the idea of uniform asymptotics in \cite{APC}, we define a number $\alpha$ such that
\begin{equation}\label{Define alpha}
 \frac 1 2   \pi i \alpha^2 =\int_{-\alpha}^\alpha(\tau^2-\alpha^2)^{1/2}d\tau=\int_{\lambda_1}^{\lambda_2}Q^{1/2}(\lambda,x)d \lambda,
\end{equation}
and a new variable $\zeta$ by
\begin{equation}\label{Define zeta}
  \int_{\alpha}^{\zeta}(\tau^2-\alpha^2)^{1/2}d\tau=\int_{\lambda_2}^\lambda Q^{1/2}(\eta,x)d\eta.
\end{equation}
Here and in \eqref{Define alpha},
the cut for the integrand on the left-hand side  is the line segment joining $-\alpha$ and  $\alpha$. The path of integration is taken   along the upper edge of the cut.
With $\alpha$ and $\zeta$ so chosen, the result   in  \cite[Theorem 1]{APC} applies. Thus we    have the following theorem.

\begin{theorem}\label{Uniform asymptotic theorem}~Given a solution $\phi(\lambda,x)$ of \eqref{sencond order equation 2},
there exist constants $c_1$ and $c_2$ such that
 \begin{equation}\label{Uniform approximate of phi}
   \Big(\frac{\zeta^2-\alpha^2}{Q(\lambda,x)}\Big)^{-\frac{1}{4}}\phi(\lambda,x)=
\Big\{[c_1+o(1)]D_{\nu}({e}^{\pi i/4}\sqrt{2x}\zeta)
+[c_2+o(1)]D_{-\nu-1}({e}^{-\pi i/4}\sqrt{2x}\zeta)\Big\}
 \end{equation}as $x\rightarrow+\infty$,  uniformly for $\lambda$ on the Stokes curves defined in
 \eqref{Stokes curves},
 where $D_{\nu}(z)$ and $D_{-\nu-1}(z)$ are   parabolic cylinder functions; cf. \cite[(12.2.5)]{OL}, with
\begin{equation}\label{Define nu}
  \nu=-\frac{1}{2}+\frac{1}{2}ix\alpha^2.
\end{equation}
\end{theorem}

%%%%%%%%%%%%%%%%%%%%%%%%%%%%%%%%%%%%%%%%%%%%%%%%%%%%%%%%%%%%%%%%%%%%%%%%%%%%%%%%%%%%%%%%%%%%%%%%%%%%%%%%%

 \section{The monodromy data for $x\rightarrow+\infty$}\label{Monod data}\setcounter{equation}{0}

We proceed  to calculate the two invariants $I_0$ and $I_1$  in  \eqref{Invariant quantity 1}-\eqref{Invariant quantity 2}
as $x\rightarrow+\infty$.
For our purpose, we need to clarify the relation  between $\zeta$ and $\lambda$ in  \eqref{Define zeta}.
\begin{lemma}\label{ralation between zeta and lambda infty}~For large $x$ and $\lambda$,
 \begin{equation} \label{ralation zeta and lambda infty}
   \frac{1}{2}ix\zeta^2-\frac{1}{2}ix\alpha^2\ln\zeta=-\frac{x}{2}\lambda+\frac{i}{4}x+\frac{x}{4}+\pi a b-\frac{1}{4}\ln x
+\frac{1}{2}\ln(b{e}^{-2is})+o(1),
 \end{equation}
 where
 \begin{equation}\label{Approximate value of alpha}
   \alpha^2=\frac{4ab+i}{x}+o\left (\frac{1}{x}\right )\quad{\rm as}~x\rightarrow\infty.
 \end{equation}
\end{lemma}

%\begin{lemma}\label{ralation between zeta and lambda orgin}~When $\lambda\rightarrow0$, for large $x$ such that $x\sqrt{\lambda}\rightarrow\infty$ and
%$x\zeta^2\rightarrow\infty$,
% \begin{equation}\label{ralation zeta and lambda orgin}
%   \frac{1}{2}ix\zeta^2-\frac{1}{2}ix\alpha^2\ln\zeta=\frac{i}{2}x\sqrt{\lambda}+\frac{i}{4}x-\frac{1}{4}\ln x-\frac{\pi i}{4}+
%\frac{1}{2}\ln(b{e}^{-2is})+o(1),
 %\end{equation}
% where $\alpha^2$ given in \eqref{Approximate value of alpha}.
%\end{lemma}

\begin{remark}
Coupling \eqref{Define nu} and \eqref{Approximate value of alpha} determines the approximate value
\begin{equation}\label{Approximate value of nu}
  \nu=2ia b-1+o(1)\quad{\rm as}~x\rightarrow\infty
\end{equation}
for the order of the parabolic cylinder function $D_{\nu}(e^{ \pi i/4}\sqrt{2x}\zeta)$ in  \eqref{Uniform approximate of phi}.
\end{remark}

\noindent{\bf Proof.}~%We will prove Lemma \ref{ralation between zeta and lambda infty}, while a similarly analysis can be applied to justify
%Lemma \ref{ralation between zeta and lambda orgin}.
 A  straightforward  integration on the left-hand side of  \eqref{Define zeta} yields
\begin{equation}\label{Left hand integral}
  \int_{\alpha}^{\zeta}(\tau^2-\alpha^2)^{1/2}d\tau=\frac{1}{2}\left \{\zeta(\zeta^2-\alpha^2)^{1/2}-\alpha^2\ln\left (\zeta+(\zeta^2-\alpha^2)^{1/2}\right )
+\alpha^2\ln\alpha\right \}.
\end{equation}
Here, the cut for the integrand is again the line segment joining $-\alpha$ and $\alpha$, and again we take  the integration path  along the upper
edge of the cut.
In view  of \eqref{Left hand integral} and   picking up the leading terms in \eqref{Define zeta},
 for large $\zeta$ we have
\begin{equation}\label{Left hand integral Asymptotic}
  \frac{1}{2}\zeta^2-\frac{1}{2}\alpha^2\ln(2\zeta)-\frac{1}{4}\alpha^2+\frac{1}{2}\alpha^2\ln(\alpha)+O\left (\alpha^4\zeta^{-2}\right )
=\int_{\lambda_2}^\lambda Q^{1/2}(\eta,x)d \eta.
\end{equation}

To calculate  the right-hand side of  \eqref{Left hand integral Asymptotic}, we split the  integration interval, so that
\begin{equation}\label{Right hand integral}
  \int_{\lambda_2}^{\lambda}Q^{1/2}(\eta,x) d \eta=
\left\{\int_{\lambda_2}^{\lambda^*}+\int_{\lambda^*}^\lambda\right\}Q^{1/2}(\eta,x) d\eta:=I_1+I_2
\end{equation}
where
\begin{equation}\label{Lambda star}
  \lambda^*=\frac{1}{2}+Tx^{-1/2},
\end{equation}
with $T$ being  a large parameter to be specified more precisely later.
                        When $\lambda$ approaches $1/2$, it follows from \eqref{Q1 function} and \eqref{Q2 function} that $Q_1\thicksim-(4ab+i)$ and $Q_2\thicksim-2(F_1(s)+F_2(s))=-4a{e}^{2is}$, where $s=\frac x 4-ab\ln\frac x 4$, $F_1$ and $F_2$ are given respectively in \eqref{Expression of F} and \eqref{Express F2}. Thus, in view of   $|{\rm Im}\{ab\}|<\frac{1}{4}$, we have
\begin{equation}\label{asymp:Q2}
 \frac{Q_2}{x^{1/2}}=o(1)\quad \mbox{as}~x\to\infty.
\end{equation}
To approximate  $I_1$,  we make the change of variables
\begin{equation}\nonumber
  \eta-\frac{1}{2}=t x^{-1/2},
\end{equation}
 replace $Q_1$ by $-(4ab+i)$,  and ignore $Q_2$. Then   for large $x$, we have
\begin{align}
I_1&=\int_{\lambda_2}^{\lambda^*}Q^{\frac{1}{2}}(\eta,x)d\eta=\frac{1}{x}\int_{\sqrt{4a b+i}}^{T}\sqrt{t^2-(4a b+i)^2}\; (1+o(1))d t\nonumber
 \\
&=\frac{T^2}{2x}-\frac{4a b+i}{4x}-\frac{4a b+i}{2x}\ln(2T)+\frac{4a b+i}{4x}\ln(4a b+i)+o\left (\frac 1 x\right ).\label{Approximate value of I1}
\end{align}
Here, as before,  the cut for the second integral  is the line segment joining $-\sqrt{4a b+i}$ and $\sqrt{4a b+i}$,
and  the path of integration is taken along the upper edge of the cut.
Taking $T=-\sqrt{4ab+i}$ in $I_1$ gives \eqref{Approximate value of alpha}.

 When $\lambda$ is large, the three terms in  the square brackets   on the righthand side of \eqref{Q2 function} are of size $o(x^{1/2})$ as $x\to\infty$. Comparing this with   \eqref{Q1 function},     one can
ignore $Q_2$ in $I_2$; cf. \eqref{sencond order equation 2} and \eqref{Right hand integral}. Accordingly  we obtain
\begin{equation}\label{Approximate value of I2}
  \begin{split}
I_2&\approx \int_{\lambda^*}^\lambda\left [-\frac{(2\eta-1)^2}{16\eta(\eta-1)}+\frac{1}{x}\left (\frac{4a b+i}{4\eta(\eta-1)}
-i\frac{2\eta-1}{8\eta(\eta-1)}\right )\frac{H'}{H}\right ]^{\frac{1}{2}}d\eta\\
&\approx i\int_{\lambda^*}^\lambda\frac{2\eta-1}{4\sqrt{\eta(\eta-1)}}\left [1-\frac{1}{ x}\frac{8}{(2\eta-1)^2}\left (\frac{4a b+i}{4}
-\frac{i(2\eta-1)}{8}\frac{H'}{H}\right )\right ]d\eta
\\
&\approx \frac{i}{2}\lambda-\frac{i}{4}+\frac{1}{4}-\frac{T^2}{2x}-\frac{i\pi(4a b+i)}{4x}-\frac{4a b+i}{4x}\ln\frac{x}{T^2}
+\frac{i}{2x}\ln\frac{ix^{1/2}}{b{e}^{-2is}}
,
\end{split}
\end{equation}with error term  being  $o(x^{-1})+O(\lambda^{-1})+O(T^4x^{-2})$, where the derivative $H'$ is taken with respect to $\eta$.
Now setting  $T<x^{1/4}$, substituting  \eqref{Approximate value of I1}  and  \eqref{Approximate value of I2}  into \eqref{Left hand integral Asymptotic},
 and combining the latter  with  \eqref{Approximate value of alpha},
we obtain Lemma \ref{ralation between zeta and lambda infty}.\hfill\qed\vskip .3cm

Now we have the relation between $\zeta$ and $\lambda$ from \eqref{Define zeta} for large $\lambda$. We also need to establish the following relation for small $\lambda$:

\begin{lemma}\label{ralation between zeta and lambda orgin}~When $\lambda\rightarrow0$, for large $x$ such that $x\sqrt{\lambda}\rightarrow\infty$ and
$x\zeta^2\rightarrow\infty$,
 \begin{equation}\label{ralation zeta and lambda orgin}
   \frac{1}{2}ix\zeta^2-\frac{1}{2}ix\alpha^2\ln\zeta=\frac{i}{2}x\sqrt{\lambda}+\frac{i}{4}x-\frac{1}{4}\ln x-\frac{\pi i}{4}+
\frac{1}{2}\ln(b{e}^{-2is})+o(1),
 \end{equation}
 where $\alpha^2$ given in \eqref{Approximate value of alpha}.
\end{lemma}

\noindent{\bf Proof.}~Now, we assume that   $\lambda\rightarrow0$   and $x$ is   large such that $x\sqrt{\lambda}\rightarrow\infty$ and $x\zeta^2\rightarrow\infty$. Let
\begin{equation}\label{lambda below star}
  \lambda_*=\frac{1}{2}-T_*x^{-1/2},
\end{equation}
where $T_*$ is a large positive constant to be specified later. We split the integral on the righthand side of \eqref{Left hand integral Asymptotic} into three parts
\begin{equation}\label{Lemma2:main:integral}
  \int_{\lambda_2}^{\lambda}Q^{1/2}(\eta,x) d \eta=
\left\{\int_{\lambda_2}^{\lambda_1}+\int_{\lambda_1}^{\lambda_*}+\int_{\lambda_*}^{\lambda}\right\}Q^{1/2}(\eta,x) d\eta:=I_0+I_1+I_2.
\end{equation}
On account of \eqref{Define alpha}, we see that
\begin{equation}\label{int:I0}
  I_0=-\frac{\pi i}{2x}\alpha^2+o\left (\frac 1 x\right ).
\end{equation}
Similar to \eqref{Approximate value of I1}, we have
\begin{align}
I_1&=\int_{\lambda_2}^{\lambda^*}Q^{\frac{1}{2}}(\eta,x)d\eta=\frac{1}{x}\int_{\sqrt{4a b+i}}^{T}\sqrt{t^2-(4a b+i)^2}\; (1+o(1))d t\nonumber
 \\
&=\frac{T_{*}^2}{2x}-\frac{4a b+i}{4x}-\frac{4a b+i}{2x}\ln(2T_*)+\frac{4a b+i}{4x}\ln(4a b+i)+o\left (\frac 1 x\right ).\label{Lm2:Approximate value of I1}
\end{align}
When $\lambda\to 0$, for large $x$, from \eqref{Q2 function} we see that
\begin{equation}\label{Lm2:Size of Q2}
  \int_{\lambda_*}^\lambda \frac{Q_2\sqrt{\lambda(\lambda-1)}}{x^{3/2}(2\lambda-1)}d\lambda
\end{equation}
is of size $o  (\sqrt{\lambda}x^{-1}  )+o  (x^{-1}\ln (\frac{x}{T_*^2}) )$.
Then, using the binomial expansion, the integral $I_2$ in \eqref{Lemma2:main:integral} is approximated as
\begin{equation}\label{Lm2:Approximate value of I2}
  \begin{split}
I_2&\approx \int_{\lambda^*}^\lambda\left [-\frac{(2\eta-1)^2}{16\eta(\eta-1)}+\frac{1}{x}\left (\frac{4a b+i}{4\eta(\eta-1)}
-i\frac{2\eta-1}{8\eta(\eta-1)}\right )\frac{H'}{H}\right ]^{\frac{1}{2}}d\eta\\
&\approx \frac{1}{2}\lambda+\frac{1}{4}-\frac{T_{*}^2}{2x}+\frac{i\pi(4a b+i)}{2x}-\frac{4a b+i}{4x}\ln\frac{x}{T_{*}^2}
+\frac{i}{2x}\ln\frac{ix^{1/2}}{b{e}^{-2is}}
,
\end{split}
\end{equation}with error term being $o(x^{-1})+o(\sqrt{\lambda}x^{-1})+o(x^{-1}\ln (\frac{x}{T_*^2}) )+O(T_*^4x^{-2})$. Now choosing  $T_*<x^{1/4}$, such that $o(x^{-1}\ln (\frac{x}{T_*^2}))=o(x^{-1})$,  substituting  \eqref{int:I0}, \eqref{Lm2:Approximate value of I1}  and  \eqref{Lm2:Approximate value of I2}  into \eqref{Lemma2:main:integral},
 and combining the latter  with  \eqref{Approximate value of alpha},
we obtain \eqref{ralation zeta and lambda orgin}, thus completing the proof of   Lemma \ref{ralation between zeta and lambda orgin}.
\hfill\qed\vskip .3cm

For large  $x$, it follows from Theorem \ref{Uniform asymptotic theorem} that
 there are two uniform  asymptotic solutions of equation  \eqref{sencond order equation 2}   $\tilde{\phi}_{\nu}$ and $\tilde{\phi}_{-\nu-1}$, namely
\begin{equation}\label{Denote of phi}
  \tilde{\phi}_{\nu}(\lambda, x)=\left(\frac{\zeta^2-\alpha^2}{Q(\lambda,x)}\right)^{\frac{1}{4}}D_{\nu}\left({e}^{\pi i/4}\sqrt{2x}\zeta\right )
\end{equation}
and
\begin{equation}\label{Denote of tilde phi}
  \tilde{\phi}_{-\nu-1}(\lambda, x)=\left (\frac{\zeta^2-\alpha^2}{Q(\lambda,x)}\right )^{\frac{1}{4}}D_{-\nu-1}\left ({e}^{-\pi i/4}\sqrt{2x}\zeta\right ),
\end{equation}
the  uniformity is  with respect to $\lambda$  on the Stokes  curves.

We denote by $\hat{\Psi}_{ij}$ the $(i,j)$ entry of $\hat{\Psi}$ and seek asymptotic solutions at infinity and at the origin.
By virtue of  \eqref{Transform psi to phi}, we have
\begin{align}
\hat{\Psi}_{11}^{(\infty)}= B^{1/2}(c_1\tilde{\phi}_{\nu}+c_2\tilde{\phi}_{-\nu-1}),\label{Asymptotic of hat Phi 11 infty}\\
\hat{\Psi}_{12}^{(\infty)}= B^{1/2}(c_3\tilde{\phi}_{\nu}+c_4\tilde{\phi}_{-\nu-1}),\label{Asymptotic of hat Phi 12 infty}
\end{align}
where $c_j$, $j=1,2,3,4$ are  constants to be determined by \eqref{Psi function expand at infty}.

Similarly, we obtain
\begin{align}
\hat{\Psi}_{11}^{(0)}= B^{1/2}(c_5\tilde{\phi}_{\nu}+c_6\tilde{\phi}_{-\nu-1}),\label{Asymptotic of hat Phi 11 orgin}\\
\hat{\Psi}_{12}^{(0)}= B^{1/2}(c_7\tilde{\phi}_{\nu}+c_8\tilde{\phi}_{-\nu-1}),\label{Asymptotic of hat Phi 12 orgin}
\end{align}
where again $c_j$, $j=5,6,7,8$ are constants to be determined.

For our purpose, we need to determine the asymptotic values of $c_j$,  $j=1,\cdots,8$. As a matter of fact, it follows from
the first equality  in  \eqref{Definition of connetion matrix E}  that
\begin{align}
\hat{\Psi}_{11}^{(\infty)}=(E_0)_{11}\hat{\Psi}_{11}^{(0)}+(E_0)_{21}\hat{\Psi}_{12}^{(0)},\label{Express hat Phi 11 infty by 0}\\
\hat{\Psi}_{12}^{(\infty)}=(E_0)_{12}\hat{\Psi}_{11}^{(0)}+(E_0)_{22}\hat{\Psi}_{12}^{(0)}.\label{Express hat Phi 12 infty by 0}
\end{align}

We obtain from  \eqref{Asymptotic of hat Phi 11 infty}, \eqref{Asymptotic of hat Phi 11 orgin}-\eqref{Asymptotic of hat Phi 12 orgin}
and \eqref{Express hat Phi 11 infty by 0}  that
\begin{equation}\nonumber
 c_1=(E_0)_{11}c_5+(E_0)_{21}c_7,~c_2=(E_0)_{11}c_6+(E_0)_{21}c_8,
\end{equation}
which in turn gives
\begin{equation}\label{Express invariant quantity 1:1}
  \frac{(E_0)_{11}}{(E_0)_{21}}=-\frac{c_1c_8-c_2c_7}{c_1c_6-c_2c_5}.
\end{equation}
Similarly calculation leads to
\begin{equation}\label{Express invariant quantity 1:2}
  \frac{(E_0)_{22}}{(E_0)_{12}}=-\frac{c_3c_6-c_4c_5}{c_3c_8-c_4c_7}.
\end{equation}
Eventually, substituting \eqref{Express invariant quantity 1:1}  and \eqref{Express invariant quantity 1:2}  into \eqref{Invariant quantity 1}  yields
\begin{equation}\label{Express invariant quantity 1}
  I_0=\frac{(c_1c_8-c_2c_7)(c_3c_6-c_4c_5)}{(c_3c_8-c_4c_7)(c_1c_6-c_2c_5)}.
\end{equation}

Now, we are in a position to calculate $c_j$ for $j=1,\cdots,8$.

From \cite[Sec.12.9]{OL}, we have the asymptotic behavior of $D_{\nu}(z)$ for $|z|\rightarrow\infty$ as follows:
\begin{equation}\label{Asymptotic behavior of D z}
  D_{\nu}(z)\sim\left\{\begin{array}{ll}
z^{\nu}{e}^{-\frac{1}{4}z^2},&\arg z\in(-\frac{3}{4}\pi,\frac{3}{4}\pi),\\
z^{\nu}{e}^{-\frac{1}{4}z^2}-\frac{\sqrt{2\pi}}{\Gamma(-\nu)}{e}^{i\pi\nu}z^{-\nu-1}{e}^{\frac{1}{4}z^2},
&\arg z\in(\frac{1}{4}\pi,\frac{5}{4}\pi).
 \end{array}
 \right.
\end{equation}
and the asymptotic behavior of $D_{-\nu-1}(iz)$ for $|z|\rightarrow\infty$  is  that
\begin{equation} \label{Asymptotic behavior of D iz}
D_{-\nu-1}(iz)\sim\left\{\begin{array}{ll}
{e}^{-\pi i(\nu+1)/2}z^{-\nu-1}{e}^{\frac{1}{4}z^2},&{\rm on}~\arg z=-\frac{1}{4}\pi,\\
-\frac{\sqrt{2\pi}}{\Gamma(\nu+1)}{e}^{-\pi i(\nu+2)/2}z^{\nu}{e}^{-\frac{1}{4}z^2},&{\rm on}~\arg z=\frac{1}{4}\pi.
 \end{array}
 \right.
\end{equation}

For $\lambda$ on the Stokes line $\arg\lambda=-\frac{\pi}{2}$ and $|\lambda|\rightarrow\infty$,
it immediately follows from  \eqref{ralation zeta and lambda infty}
that $\zeta^2\thicksim i\lambda$, if we take $i={e}^{\frac{\pi i}{2}}$, then we have $\arg\zeta\sim0$. Therefore,
$\arg({e}^{\pi i/4}\sqrt{2x}\zeta)\thicksim\frac{\pi}{4}$ and $\arg({e}^{-\pi i/4}\sqrt{2x}\zeta)\thicksim-\frac{\pi}{4}$ for $x>0$. Taking
$-1={e}^{\pi i}$, from  \eqref{sencond order equation 2}  we have $Q^{-1/4}\thicksim\sqrt{2}{e}^{-\pi i/4}$ as $|\lambda|\rightarrow\infty$. Since
$(\zeta^2-\alpha^2)^{1/4}\thicksim\zeta^{1/2}$ as $|\lambda|\rightarrow\infty$,   using the appropriate asymptotic formulas of $D_{\nu}(z)$ in
 \eqref{Asymptotic behavior of D z}, and in view of  \eqref{Denote of phi}  and  \eqref{ralation zeta and lambda infty}, we have
\begin{equation}
 \tilde{\phi}_{\nu}\thicksim A_0{e}^{\frac{x}{2}\lambda},~{\rm as}~|\lambda|\rightarrow\infty.
\label{Asymptotic of phi}
\end{equation}
Here use has been  made of $s=\frac{x}{4}-ab\ln\frac{x}{4}$, and
\begin{equation}\nonumber
  A_0=\sqrt{2}{e}^{-\pi i/4}x^{-1/4}{e}^{-\frac{1}{4}x}2^{3iab-1/2}b^{-1/2}
{e}^{-\pi ab}{e}^{\frac{\pi i}{4}\nu}.
\end{equation}
Similarly, by applying the appropriate asymptotic formulas of $D_{\nu}(z)$ in
 \eqref{Asymptotic behavior of D z}, we obtain from \eqref{Denote of tilde phi}  and  \eqref{ralation zeta and lambda infty} that
\begin{equation}\label{Asymptotic of tilde phi}
  \tilde{\phi}_{-\nu-1}\thicksim B_0{e}^{-\frac{x}{2}\lambda},~{\rm as}~|\lambda|\rightarrow\infty,
\end{equation}
where
\begin{equation}\nonumber
  B_0=\sqrt{2}{e}^{-\pi i/4}x^{-1/4}{e}^{\frac{1}{4}x}2^{-3iab}b^{1/2}
{e}^{\pi ab}{e}^{\frac{\pi i}{4}(\nu+1)}.
\end{equation}

On the other hand, by virtue of \eqref{Gauge transformation}, we obtain from \eqref{Psi function expand at infty} that
\begin{align}
\hat{\Psi}_{11}^{(\infty)}&\thicksim\frac{1}{\sqrt{2}}u^{-\frac{1}{2}}{e}^{\frac{x}{2}\lambda},\label{Asymptoic hat Psi:11}\\
\hat{\Psi}_{12}^{(\infty)}&\thicksim\frac{1}{\sqrt{2}}iu^{\frac{1}{2}}{e}^{-\frac{x}{2}\lambda},\label{Asymptoic hat Psi:12}
\end{align}
Thus, applying \eqref{Asymptotic of phi}  and  \eqref{Asymptotic of tilde phi}, a combination of   \eqref{Asymptotic of hat Phi 11 infty},  \eqref{Asymptotic of hat Phi 12 infty},  \eqref{Asymptoic hat Psi:11}  and  \eqref{Asymptoic hat Psi:12}  implies
\begin{equation}\label{Value of c2 c3}
  c_2=c_3=0.
\end{equation}
Substituting  \eqref{Value of c2 c3}  into  \eqref{Express invariant quantity 1}  gives the simplified version
\begin{equation}\label{Express invariant quantity 2}
  I_0=\frac{c_5c_8}{c_6c_7}.
\end{equation}

Our next task is to determine $c_j$ for $j=5,6,7,8$.

By using the result from \cite[Prop.7.1]{FA4}, we have
\begin{equation}\label{Hat Phi at lambda=0}
  \hat{\Psi}^{(0)}\thicksim M_1
\left(\begin{array}{cc}
v+\Theta &~v+\Theta\\
v &~v+\Theta
\end{array}\right)(I+o(1))\; \Phi\left (\frac{1}{2}x\sqrt{\lambda},\Theta\right )\; u^{-\frac{1}{2}\sigma_3},
\end{equation}
as $\lambda\rightarrow0$ and  $x\rightarrow\infty$,
where $M_1$ is defined in  \eqref{M1}  and
\begin{equation}\label{Model Bessel function}
  \Phi(x,\Theta)=x\left(\begin{array}{cc}
J_{\Theta-1}(x) &~J_{1-\Theta}(x)\\
J_{\Theta+1}(x) &~J_{-(\Theta+1)}(x)
\end{array}\right),
\end{equation}
where $J_{\mu}(x)$  is the Bessel function of the first kind,  and $\Theta\notin\mathbb{Z}$, as mentioned earlier.

For $\arg\lambda=2\pi$, choosing $\lambda\rightarrow0$ such that $x\sqrt{\lambda}\rightarrow\infty$ and $x\zeta^2\rightarrow\infty$,
by virtue of the asymptotic expansion of $v(x)$ in \eqref{Expanding of v at infty}   and the behavior for Bessel function
\begin{equation}\nonumber
  J_{\mu}(z)=\sqrt{\frac{2}{\pi z}}\left [\cos\left (z-\frac{1}{2}\mu\pi-\frac{1}{4}\pi\right )
  -\frac{1}{2z}\left(\mu^2-\frac{1}{4}\right)\sin\left(z-\frac{1}{2}\mu\pi-\frac{1}{4}\pi\right)
+O\left (z^{-2}\right )\right ],
\end{equation}
as $z\rightarrow\infty$ for $|\arg z|<\pi$; cf. \cite{OL, Whittaker and Watson},   a straightforward calculation enables us to
write \eqref{Hat Phi at lambda=0}  as the following form
\begin{equation}\nonumber
 \hat{\Psi}^{(0)}\thicksim \frac{-\Theta}{2\sqrt{\pi}}x^{1/2}\lambda^{-1/4} M_1
\left(\begin{array}{cc}
\cos(\frac{1}{2}x\sqrt{\lambda}-\frac{1}{2}\Theta\pi-\frac{1}{4}\pi)&-\cos(\frac{1}{2}x\sqrt{\lambda}+\frac{1}{2}\Theta\pi-\frac{1}{4}\pi)\\
\cos(\frac{1}{2}x\sqrt{\lambda}-\frac{1}{2}\Theta\pi-\frac{1}{4}\pi)&-\cos(\frac{1}{2}x\sqrt{\lambda}+\frac{1}{2}\Theta\pi-\frac{1}{4}\pi)
\end{array}\right)u^{-\frac{1}{2}\sigma_3},
\end{equation}
which gives us
\begin{equation}
 \hat{\Psi}_{11}^{(0)}\thicksim-\frac{\Theta}{4\sqrt{\pi}}x^{1/2}u^{-\frac{1}{2}}\lambda^{-1/4}(1+i)
 \left (e^{-\frac{i}{2}\Theta\pi-\frac{i}{4}\pi}e^{\frac{i}{2}x\sqrt{\lambda}}
+e^{\frac{i}{2}\Theta\pi+\frac{i}{4}\pi}e^{-\frac{i}{2}x\sqrt{\lambda}}\right )
\label{asymptotic hat Phi 11 0}
\end{equation}
and
\begin{equation}\label{asymptotic hat Phi 12 0}
  \hat{\Psi}_{12}^{(0)}\thicksim\frac{\Theta}{4\sqrt{\pi}}x^{1/2}u^{\frac{1}{2}}\lambda^{-1/4}(1+i)
  \left (e^{\frac{i}{2}\Theta\pi-\frac{i}{4}\pi}e^{\frac{i}{2}x\sqrt{\lambda}}
+e^{-\frac{i}{2}\Theta\pi+\frac{i}{4}\pi}e^{-\frac{i}{2}x\sqrt{\lambda}}\right )
\end{equation}
as $x\sqrt{\lambda}\rightarrow\infty$.

On the other hand, for $\arg\lambda=2\pi$, choosing $\lambda\rightarrow0$ such that $x\sqrt{\lambda}\rightarrow\infty$ and $x\zeta^2\rightarrow\infty$,
then we have from \eqref{ralation zeta and lambda orgin}  that
$\zeta^2\thicksim\sqrt{\lambda}$, which gives us $\arg\zeta\thicksim\frac{1}{2}\pi$. Hence,
$\arg({e}^{\pi i/4}\sqrt{2x}\zeta)\thicksim\frac{3\pi}{4}$ and $\arg({e}^{-\pi i/4}\sqrt{2x}\zeta)\thicksim\frac{\pi}{4}$ for $x>0$.
From  \eqref{sencond order equation 2}  we have $Q^{-1/4}\thicksim2\lambda^{1/4}$ as $|\lambda|\rightarrow0$.
By using the appropriate asymptotic formulas of $D_{\nu}(z)$ in
 \eqref{Asymptotic behavior of D z}, we obtain from  \eqref{Denote of phi}  and  \eqref{ralation zeta and lambda orgin} that
\begin{equation}\label{Asymptotic of phi 0}
 \tilde{\phi}_{\nu}\thicksim2\lambda^{1/4}x^{-1/4}(C_0{e}^{-\frac{i}{2}x\sqrt{\lambda}}-D_0{e}^{\frac{i}{2}x\sqrt{\lambda}}),
\end{equation}
where
\begin{equation}\label{coefficient in psi as lambda 0}
  C_0={e}^{\pi i(\nu+1)/4}2^{3iab-1/2}b^{-1/2},~D_0=\frac{\sqrt{2\pi}}{\Gamma(-\nu)}{e}^{\frac{3i\pi}{4}\nu}
{e}^{-\frac{\pi i}{2}}2^{-3iab}b^{1/2}.
\end{equation}

Similarly, by further applying the appropriate asymptotic formulas for $D_{-\nu-1}(-iz)$ in
 \eqref{Asymptotic behavior of D iz}, we obtain from  \eqref{Denote of tilde phi}  and  \eqref{ralation zeta and lambda orgin} that
\begin{equation}\label{Asymptotic of tilde phi 0}
  \tilde{\phi}_{-\nu-1}\thicksim2\lambda^{1/4}x^{-1/4}F_0{e}^{\frac{i}{2}x\sqrt{\lambda}}
\end{equation}
with
\begin{equation}\label{coefficient in tilde psi 0}
  F_0={e}^{\frac{\pi i}{4}\nu}2^{-3iab}b^{1/2}.
\end{equation}

Substituting  \eqref{Asymptotic of phi 0}  and  \eqref{Asymptotic of tilde phi 0}  into  \eqref{Asymptotic of hat Phi 11 orgin}  and
using  \eqref{asymptotic hat Phi 11 0}, we obtain
\begin{align}
&B^{1/2}2\lambda^{1/4}x^{-1/4}c_5 \thicksim-\frac{\Theta}{4\sqrt{\pi}}x^{1/2}u^{-\frac{1}{2}}\lambda^{-1/4}(1+i)
{e}^{\frac{i}{2}\Theta\pi+\frac{i}{4}\pi}C_0^{-1},\label{c5}\\
&B^{1/2}2\lambda^{1/4}x^{-1/4}c_6F_0 \thicksim-\frac{\Theta}{4\sqrt{\pi}}x^{1/2}u^{-\frac{1}{2}}\lambda^{-1/4}(1+i)
({e}^{-\frac{i}{2}\Theta\pi-\frac{i}{4}\pi}+{e}^{\frac{i}{2}\Theta\pi+\frac{i}{4}\pi}D_0C_0^{-1}).\label{c6}
\end{align}
Similarly, substituting \eqref{Asymptotic of phi 0}  and  \eqref{Asymptotic of tilde phi 0} into  \eqref{Asymptotic of hat Phi 12 orgin}  and
using  \eqref{asymptotic hat Phi 12 0}, we obtain
\begin{align}
& B^{1/2}2\lambda^{1/4}x^{-1/4}c_7  \thicksim\frac{\Theta}{4\sqrt{\pi}}x^{1/2}u^{\frac{1}{2}}\lambda^{-1/4}(1+i)
{e}^{-\frac{i}{2}\Theta\pi+\frac{i}{4}\pi}C_0^{-1},\label{c7}\\
& B^{1/2}2\lambda^{1/4}x^{-1/4}c_8F_0 \thicksim\frac{\Theta}{4\sqrt{\pi}}x^{1/2}u^{\frac{1}{2}}\lambda^{-1/4}(1+i)
({e}^{\frac{i}{2}\Theta\pi-\frac{i}{4}\pi}+{e}^{-\frac{i}{2}\Theta\pi+\frac{i}{4}\pi}D_0C_0^{-1}).\label{c8}
\end{align}
Therefore, taking ratios    from  \eqref{c5}  and  \eqref{c7}, we have
\begin{equation}\label{c5:c7}
  \frac{c_5}{c_7}=-{e}^{\pi i\Theta}u^{-1},
\end{equation}
and a parallel result follows from \eqref{c6}  and  \eqref{c8} reads
\begin{equation}\label{c8:c6}
  \frac{c_8}{c_6}=-u\frac{e^{\frac{i}{2}\Theta\pi-\frac{i}{4}\pi}+e^{-\frac{i}{2}\Theta\pi+\frac{i}{4}\pi}D_0C_0^{-1}}{
{e^{-\frac{i}{2}\Theta\pi-\frac{i}{4}\pi}+e^{\frac{i}{2}\Theta\pi+\frac{i}{4}\pi}D_0C_0^{-1}}},
\end{equation}
where
\begin{equation}\label{D0 C0^(-1)}
  D_0C_0^{-1}=-\frac{2\sqrt{\pi}b}{\Gamma(1-2iab)}e^{-\frac{i}{4}\pi}e^{-\pi ab}2^{-6iab}.
\end{equation}
Substituting  \eqref{c5:c7}  and  \eqref{c8:c6}  into  \eqref{Express invariant quantity 2}, we obtain the first invariant quantity for large $x$
\begin{equation}\label{Invariant quantity 1 large x}
  I_0=\frac{i{e}^{\pi i\Theta}-D_0C_0^{-1}}{i{e}^{-\pi i\Theta}-D_0C_0^{-1}}.
\end{equation}

For large $x$, in the same manner, the other invariant can also be calculated as
\begin{equation}\label{Invariant quantity 2 large x}
  I_1=\frac{i{e}^{\pi i\Theta}-\frac{2\sqrt{\pi}a}{\Gamma(1+2iab)}{e}^{\frac{i}{4}\pi}{e}^{-\pi ab}2^{6iab}}
{i{e}^{-\pi i\Theta}-\frac{2\sqrt{\pi}a}{\Gamma(1+2iab)}{e}^{\frac{i}{4}\pi}{e}^{-\pi ab}2^{6iab}}.
\end{equation}

A combination  of \eqref{Invariant quantity 1 small x}  and  \eqref{Invariant quantity 1 large x}  with  \eqref{D0 C0^(-1)}  gives  \eqref{Connection formula b},
and a combination of \eqref{Invariant quantity 2 small x}  and  \eqref{Invariant quantity 2 large x}  gives \eqref{Connection formula a}.

%%%%%%%%%%%%%%%%%%%%%%%%%%%%%%%%%%%%%%%%%%%%%%%%%%%%%%%%%%%%%%%%%%%%%%%%%%%%%%%%%%%%%%%%%

\section*{Acknowledgements}
  Yu-Qiu Zhao  was supported in part by the National
Natural Science Foundation of China under grant number
10871212.

%%%%%%%%%%%%%%%%%%%%%%%%%%%%%%%%%%%%%%%%%%%%%%%%%%%%%%%%%%%%%%%%%%%%%%%%%%%%%%%%%%%%%%%%%%

\end{document}